\documentclass[11pt,reqno]{amsart}
\usepackage{mathrsfs}
\usepackage{amsfonts}
\usepackage{amssymb,amsmath,amscd,amsthm,amsfonts}
\usepackage[colorlinks=true]{hyperref}
\newcommand{\bqa}{\begin{equation}}
\newcommand{\eqa}{\end{equation}}
\newcommand{\bea}{\begin{eqnarray}}
\newcommand{\eea}{\end{eqnarray}}
\newcommand{\bna}{\begin{eqnarray*}}
\newcommand{\ena}{\end{eqnarray*}}

\def\bz{{\mathbb Z}}
\def\br{{\mathbb R}}

\def\sl2z{SL(2,\bz)}
\def\psl2z{PSL(2,\bz)}
\def\gl2r{GL(2,\br)} 

\theoremstyle{definition}

\begin{document}
\noindent {\Large \bf Weighted average values of automorphic $L$-functions
}

\bigskip

\noindent {\bf Wei Liu \footnote{ This work is supported
by the National Natural Science Foundation of China (Grant No.11871306 ).}}

\bigskip

\noindent{\footnotesize \bf Abstract} \,{\footnotesize
Let $S_2^*(q)$ be the set of primitive Hecke eigenforms of weight 2 and prime level $q$.
For $p$ prime and $t\in \mathbb{R}$, we prove asymptotic formulas for the sums
\bna
\mathcal {A}(p^j,q,t)=\sum_{f\in S_2^*(q)} L\left(\frac{1}{2}+it,f\right)^2\lambda_f(p^j),\qquad j=1,2,
\ena
where $\lambda_f(p^j)$ is the $p^j$-th normalized Fourier coefficient of $f$
and $L(s,f)$ is the $L$-function associated to $f$.
}

\medskip



\bigskip

\noindent {\bf 1 \quad Introduction}

\setcounter{section}{1} \setcounter{equation}{0}

\medskip

Critical $L$-values have been studied intensively in number theory.
One is particularly interested in establishing asymptotic formulas
for moments of various families of $L$-functions on the critical line.
Most of the time asymptotic formulas for harmonic moments over families of cusp forms
(sums weighted by the inverse of the Petersson inner product $\langle f, f\rangle$)
are sufficient for applications such as non-vanishing problems of $L$-functions, the sub-convexity problem, the problem of
non-existence of Landau-Siegel zeros and so on (see Duke \cite{D}, Kowalski and Michel \cite{KM1}, Conrey and Iwaniec \cite{CI}, Iwaniec and Sarnak \cite{IS}, Luo \cite{L} and Balkanova and Frolenkov \cite{BF} ). While in a few applications the
non-harmonic moments are needed (see VanderKam \cite{V} for example).
In this paper, we study the following
weighted average values of automorphic $L$-functions on the critical line
\bna
\mathcal {A}(p^j,q,t)=\sum_{f\in S_2^*(q)} L\left(\frac{1}{2}+it,f\right)^2\lambda_f(p^j), \qquad j=1,2,
\ena
where $S_2^*(q)$ denotes the set of primitive Hecke eigenforms of weight 2 and prime level $q$,
$\lambda_f(p^j)$ is the $p^j$-th normalized Fourier coefficient of $f$ and $L(s,f)$ is the $L$-function associated to $f$.

The main results of this paper are the following theorems.

\medskip

\noindent{\bf Theorem 1.1} \, {\it Let $p$ and $q$ be prime. Then for any $\epsilon>0$ and $p<q^{\frac{1}{55}-\epsilon}$,
we have
\bna
\mathcal {A}(p,q,0)&=&\left(1+\frac{1}{p}\right)\frac{1-q^{-1}}{1+p^{-2}}
\frac{\zeta^3(2)}{6\zeta(4)}\frac{q}{\sqrt{p}}
\left(\log\frac{q}{4\pi^2p}
+\frac{2\log{p}}{1+p^2}
+6\frac{\zeta'(2)}{\zeta(2)}-4\frac{\zeta'(4)}{\zeta(4)}\right)
+O_\epsilon\left(p^{\frac{39}{32}}q^{\frac{31}{32}+\epsilon}\right).
\ena
Given $0\neq t\in \mathbb{R}$. Then for any $\epsilon>0$ and $p<q^{\frac{1}{55}-\epsilon}$, we have
\bna
\mathcal {A}(p,q,t)&=&
\frac{\Gamma(1-it)^2}{\Gamma(1+it)^2}\left(\frac{q}{4\pi^2}\right)^{-2it}
\frac{\zeta_q(1-2it)}{6(1+p^{-2+2it})}\frac{\zeta^3(2-2it)}{\zeta(4-4it)}\frac{q(p+1)}{p^{\frac{3}{2}-it}}\nonumber\\
&&+\frac{\zeta_q(1+2it)}{6(1+p^{-2-2it})}\frac{\zeta^3(2+2it)}{\zeta(4+4it)}\frac{q(p+1)}{p^{\frac{3}{2}+it}}
+O_{t,\epsilon}\left(p^{\frac{39}{32}}q^{\frac{31}{32}+\epsilon}\right),
\ena
where $\zeta_q(s)=(1-q^{-s})\zeta(s)$.
}

\medskip

\noindent{\bf Theorem 1.2} \, {\it Let $p$ and $q$ be prime.
Then for any $\epsilon>0$ and $p<q^{\frac{1}{110}-\epsilon}$,
we have
\bna
\mathcal {A}(p^2,q,0)&=&
\frac{(q-1)\zeta^3(2)}{12p^2\zeta(4)}
\left(p+\frac{(p+1)(3-p^{-2})}{1+p^{-2}}\right)\left(\log\frac{q}{4\pi^2p^2}
+6\frac{\zeta'(2)}{\zeta(2)}
-4\frac{\zeta'(4)}{\zeta(4)}\right)\\
&+&\frac{(q-1)\zeta^3(2)\log p}{6p\zeta(4)}
\left(1+\frac{4p(p+1)(3-p^{-2})}{(1+p^{-2})(p^2+1)(3p^2-1)}
\right)\\
&+&O_\epsilon(p^{\frac{39}{16}+\epsilon}q^{\frac{31}{32}+\epsilon}).
\ena
Given $0\neq t\in \mathbb{R}$. Then for any $\epsilon>0$, we have
\bna
\mathcal {A}(p^2,q,t)&=&
\frac{\zeta_q(1+2it)\zeta^3(2+2it)}{12\zeta(4+4it)}
\left(\frac{q}{p}+\frac{q(p+1)(3-p^{-2-2it})}{p^{2+2it}+1}\right)\\
&+&\frac{\Gamma(1-it)^2}{\Gamma(1+it)^2}
\left(\frac{q}{4\pi^2}\right)^{-2it}
\frac{\zeta_q(1-2it)\zeta^3(2-2it)}{12\zeta(4-4it)}
\left(\frac{q}{p}+\frac{q(p+1)(3-p^{-2+2it})}{p^{2-2it}+1}\right)\\
&+&O_{t,\epsilon}\left(p^{\frac{39}{16}+\epsilon}q^{\frac{31}{32}+\epsilon}\right),
\ena
where $\zeta_q(s)=(1-q^{-s})\zeta(s)$.
}
\medskip

To prove the above results, we follow closely VanderKam \cite{V}. The bound $p<q^{1/55-\epsilon}$
comes from the exponent $39/32$ in Lemma 2.3 due to VanderKam \cite{V}, which in turn
comes from Burgess's bound for Dirichlet $L$-functions, that is, $L(s,\chi_D)\ll |s|^A D^{3/16+\epsilon}$
for some $A>0$ and any $\epsilon>0$. Under Lindel\"{o}f hypothesis for Dirichlet $L$-functions,
the error term $p^{39/32}q^{31/32+\epsilon}$ can be improved to $p^{9/8}q^{7/8+\epsilon}$ and
the range of $p$ can be extended to $p<q^{\frac{1}{13}-\epsilon}$ (see the discussions after Lemma 2.3).

By  using the methods in Kowalski and Michel \cite{KM2}, Bui \cite{B} proved an asymptotic formula for the second harmonic moment
\bna
\mathcal {H}(p,q)=\sum_{f\in S_2^*(q)} \omega_f^{-1}L\left(\frac{1}{2},f\right)^2\lambda_f(p),
\ena
where $\omega_f=4\pi \langle f, f\rangle$. Moreover,
he showed that for any $\epsilon>0$ and $p<q^{1/2-\epsilon}$,
\bna
\mathcal {H}(p,q)-\sqrt{\frac{p}{q}}\mathcal {H}(q,p)=\frac{2}{\sqrt{p}}\log\frac{q}{4\pi^2p}
+O\left(p^{\frac{1}{2}+\epsilon}q^{-\frac{1}{2}}\right).
\ena
From the proof of Theorem 1.1, we can see that
\bna
\sum_{f\in S_2^*(q)} L\left(\frac{1}{2},f\right)^2=
\frac{\zeta^3(2)}{12\zeta(4)}q\log \frac{q}{4\pi^2}+
\frac{\zeta^3(2)}{6\zeta(4)}\left(3\frac{\zeta'(2)}{\zeta(2)}-
2\frac{\zeta'(4)}{\zeta(4)}\right)q
+O_{\epsilon}\left(q^{\frac{31}{32}+\epsilon}\right).
\ena
Thus for $p, q$ prime,
\bna
|\mathcal {A}(q,p,0)|\leq \tau (q)\sum_{f\in S_2^*(p)} L\left(\frac{1}{2},f\right)^2 \ll p\log p.
\ena
Then by Theorem 1.1, for any $\epsilon>0$ and $p<q^{\frac{1}{55}-\epsilon}$,
we have
\bna
\mathcal {A}(p,q,0)-\sqrt{\frac{p}{q}}\mathcal {A}(q,p,0)=\frac{\zeta^3(2)}{6\zeta(4)}\frac{q}{\sqrt{p}}
\log\frac{q}{4\pi^2p}+\frac{\zeta^3(2)}{3\zeta(4)}\left(3\frac{\zeta'(2)}{\zeta(2)}-
2\frac{\zeta'(4)}{\zeta(4)}\right)\frac{q}{\sqrt{p}}
+O_{\epsilon}\left(p^{\frac{39}{32}}q^{\frac{31}{32}+\epsilon}\right).
\ena

\medskip

\noindent {\bf 2 \quad Preliminaries}

\setcounter{section}{2} \setcounter{equation}{0}

\medskip

Let $S_2^*(q)$ be the set of primitive Hecke eigenforms of weight 2 and prime level $q$.
Any $f\in S_2^*(q)$ has the Fourier expansion
\bna
f(z)=\sum_{n\geq 1}\sqrt{n}\lambda_f(n)e(nz), \quad \lambda_f(1)=1.
\ena
The Fourier coefficients $\lambda_f(n)$ satisfy
\bea
\lambda_f(m)\lambda_f(n)=\sum_{d|(m,n) \atop (d,q)=1}\lambda_f\left(mn/d^2\right).
\eea
For $\Re(s)>1$, the $L$-function associated to $f$ is defined by
\bna
L(s,f)=\sum_{n\geq 1}\lambda_f(n)n^{-s}
\ena
which has an analytic continuation to all $s\in \mathbb{C}$ and satisfies
the functional equation
\bna
\Lambda(s,f)=\left(\frac{\sqrt{q}}{2\pi}\right)^{s}\Gamma\left(s+\frac{1}{2}\right)L(s,f)=\varepsilon_f\Lambda(1-s,f),
\qquad \varepsilon_f=\pm 1.
\ena

Note that for $\Re(s)>1$,
\bna
L(s,f)^2=\sum_{d\geq 1 \atop (d,q)=1}\frac{1}{d^{2s}}\sum_{n\geq 1}\frac{\tau(n)\lambda_f(n)}{n^{s}},
\ena
where $\tau(n)$ is the divisor function.
By Iwaniec and Kowalski \cite{IK} (see Theorem 5.3), $L\left(1/2+it,f\right)^2$
has the following approximate functional equation.

\medskip

\noindent {\bf Lemma 2.1} \, {\it For any $t\in \mathbb{R}$, we have
\bea \label{L2}
L\left(\frac{1}{2}+it,f\right)^2&=&\Gamma(1+it)^{-2}\sum_{d\geq 1 \atop (d,q)=1}\frac{1}{d^{1+2it}}
\sum_{n\geq 1}\frac{\tau(n)\lambda_f(n)}{n^{\frac{1}{2}+it}}W_t\left(\frac{4\pi^2 nd^2}{q}\right)\nonumber\\
&&+\Gamma(1+it)^{-2}\left(\frac{q}{4\pi^2}\right)^{-2it}\sum_{d\geq 1 \atop (d,q)=1}\frac{1}{d^{1-2it}}
\sum_{n\geq 1}\frac{\tau(n)\lambda_f(n)}{n^{\frac{1}{2}-it}}W_{-t}\left(\frac{4\pi^2 nd^2}{q}\right),
\eea
where
\bea \label{WtY}
W_t(Y)=\frac{1}{2\pi i}\int\limits_{(2)}
Y^{-u}\Gamma(1+it+u)^2e^{u^2}\frac{\mbox{d}u}{u}.
\eea
}

\medskip

\noindent{\bf Lemma 2.2}\, {\it For $0<Y\leq 1$, we have
\bna
W_t(Y)=\Gamma(1+it)^2+O_\epsilon(Y^{1-\epsilon})
\ena
for any $\epsilon>0$. For $Y>1$, we have
\bna
W_t(Y)\ll e^{-C\sqrt{Y}}
\ena
for any constant $C>0$.}

\begin{proof}
The proof of this lemma is similar as Proposition 5.4 in \cite{IK}.
For $0<Y\leq 1$, we shift the contour of integration in (2.3) to $\Re(u)=-1+\epsilon$, passing a simple pole at $u=0$, and get
\bna
W_t(Y)&=&\Gamma(1+it)^2+
\frac{1}{2\pi i}\int\limits_{(-1+\epsilon)}
Y^{-u}
\Gamma(1+it+u)^2e^{u^2}\frac{\mbox{d}u}{u}\\
&=&\Gamma(1+it)^2+
\frac{1}{2\pi}\int\limits_{-\infty}^{\infty}
Y^{1-\epsilon-iv}
\Gamma(\epsilon+i(t+v))^2e^{(-1+\epsilon+iv)^2}\frac{\mbox{d}v}{-1+\epsilon+iv}\\
&=&\Gamma(1+it)^2+
O_\epsilon\left(Y^{1-\epsilon}\int\limits_{-\infty}^{\infty}
\left|\Gamma(\epsilon+i(t+v))\right|^2
e^{-v^2}\frac{\mbox{d}v}{1+|v|}\right)\\
&=&\Gamma(1+it)^2+O_\epsilon(Y^{1-\epsilon})
\ena
for any $\epsilon>0$.
For $Y>1$, we move the line of integration to $\Re(u)=Y$, then we  can obtain
\bna
W_t(Y)&=&\frac{1}{2\pi i}\int\limits_{(Y)}
Y^{-u}
\Gamma(1+it+u)^2e^{u^2}\frac{\mbox{d}u}{u}\\
&\ll&Y^{-Y}\Gamma(Y+1)\\
&\ll&Y^{-Y}Y^{Y+\frac{1}{2}}e^{-Y}\\
&\ll&e^{-C\sqrt Y}
\ena
for any $C>0$.
\end{proof}

\medskip

Define $Tr(T_n)=\sum\limits_{f\in S_2^*(q)}\lambda_f(n)$. The following result due to VanderKam \cite{V}
(see Lemma 6.1) is crucial to the proof of Theorem 1.1.

\medskip

\noindent {\bf Lemma 2.3} \, {\it Let $a_n$ be any sequence of numbers and $q$ be prime.
Given any $c>1$ and any $B<q^2/4c$,
\bna
\sum_{n=B}^{cB}a_nTr(T_n)=\frac{q}{12}\sum_{n=B}^{cB}\frac{a_n}{\sqrt{n}}\delta_{n=\diamond}+
O_{c,\epsilon}\left(a_{max}\left(B^{\frac{7}{4}}q^{-\frac{1}{2}}+B^{\frac{39}{32}+\epsilon}q^{\frac{1}{4}}\right)\right),
\ena
where $k=\diamond$ indicates that $k$ is a square and $a_{max}$ is the largest value of $|a_n|$ in the given range.
}

\medskip

As remarked by VanderKam \cite{V}, under Lindel\"{o}f hypothesis for Dirichlet $L$-functions, the exponent $39/32$ can be improved to $9/8$.

\medskip

\noindent {\bf 3 \quad Proof of Theorem 1.1}

\setcounter{section}{3} \setcounter{equation}{0}

\medskip
We follow closely VanderKam \cite{V}.
Note that $\overline{W_t(Y)}=W_{-t}(Y)$.
By the approximate functional equation in (2.2), we have
\bea
\mathcal {A}(p,q,t) \label{Apqt}
=\Gamma(1+it)^{-2}\mathcal {B}(p,q,t)+\Gamma(1+it)^{-2}\left(\frac{q}{4\pi^2}\right)^{-2it}\overline{\mathcal {B}(p,q,t)},
\eea
where
\bna
\mathcal {B}(p,q,t)=\sum_{d\geq 1 \atop (d,q)=1}\frac{1}{d^{1+2it}}
\sum_{n\geq 1}\frac{\tau(n)}{n^{\frac{1}{2}+it}}W_t\left(\frac{4\pi^2 nd^2}{q}\right)
\sum_{f\in \mathcal {S}_2^*(q)}\lambda_f(n)\lambda_f(p).
\ena
Throughout the proof, we assume that $p<q^{1/3}$. Then by (2.1),
\bna
\lambda_f(n)\lambda_f(p)=\sum_{d|(n,p) \atop (d,q)=1}\lambda_f\left(\frac{np}{d^2}\right)
=\left\{\begin{array}{ll} \lambda_f(np), & \mbox {if $p{\not|} n$},\\
\lambda_f(np)+\lambda_f\left(n/p\right), &\mbox {if $p|n$.}\end{array}\right.
\ena
It follows that
\bea \label{Bpqt}
\mathcal {B}(p,q,t)=\mathcal {B}_1(p,q,t)+\mathcal {B}_2(p,q,t),
\eea
where
\bea
\mathcal {B}_1(p,q,t)&=&\sum_{d\geq 1 \atop (d,q)=1}\frac{1}{d^{1+2it}}
\sum_{(n,p)=1}\frac{\tau(n)}{n^{\frac{1}{2}+it}}W_t\left(\frac{4\pi^2 nd^2}{q}\right)
\sum_{f\in \mathcal {S}_2^*(q)}\lambda_f(np),\\
\mathcal {B}_2(p,q,t)&=&\sum_{d\geq 1 \atop (d,q)=1}\frac{1}{d^{1+2it}}
\sum_{p|n}\frac{\tau(n)}{n^{\frac{1}{2}+it}}W_t\left(\frac{4\pi^2 nd^2}{q}\right)
\sum_{f\in \mathcal {S}_2^*(q)}\left(\lambda_f(np)+\lambda_f\left(\frac{n}{p}\right)\right).
\eea

\medskip

{\bf Evaluation of $\mathcal {B}_1(p,q,t)$}. We first compute $\mathcal {B}_1(p,q,t)$ in (3.3).
Since $\sum_{d|m}\mu(d)=1$ if $m=1$ and equals 0 otherwise, we have
\bna
\mathcal {B}_1(p,q,t)&=&\sum_{d\geq 1 \atop (d,q)=1}\frac{1}{d^{1+2it}}
\sum_{n\geq 1}\frac{\tau(n)}{n^{\frac{1}{2}+it}}W_t\left(\frac{4\pi^2 nd^2}{q}\right)\sum_{d|(n,p)}\mu(d)
\sum_{f\in \mathcal {S}_2^*(q)}\lambda_f(np)\\
&=&p^{\frac{1}{2}+it}\sum_{d\geq 1 \atop (d,q)=1}\frac{1}{d^{1+2it}}\sum_{k \geq 1}\frac{\tau(k/p)}{k^{\frac{1}{2}+it}}
W_t\left(\frac{4\pi^2 kd^2}{pq}\right)\sum_{d|(k/p,p)}\mu(d) Tr(T_k).
\ena
By Lemma 2.2, we can assume that $kd^2\leq \exp([\log pq^{11/10}])\leq q^2$,
where $[x]$ denotes the greatest integer not exceeding $x$.
Then we can apply Lemma 2.3 to obtain
\bna
\mathcal {B}_1(p,q,t)&=&\frac{qp^{\frac{1}{2}+it}}{12}\sum_{d\geq 1 \atop (d,q)=1}\frac{1}{d^{1+2it}}
\sum_{k =\diamond}\frac{\tau(k/p)}{k^{1+it}}W_t\left(\frac{4\pi^2 kd^2}{pq}\right)\sum_{d|(k/p,p)}\mu(d)
\nonumber\\&&+
O\left(p^{\frac{1}{2}}\sum_{B}
\left(B^{\frac{7}{4}}q^{-\frac{1}{2}}+B^{\frac{39}{32}+\epsilon}q^{\frac{1}{4}}\right)
\sum_{d\geq 1 \atop (d,q)=1}\frac{1}{d}
\max_{B\leq k\leq eB}\frac{\tau(k/p)}{\sqrt{k}}\left|
W_t\left(\frac{4\pi^2 kd^2}{pq}\right)\right|\right),
\ena
where the summation over $B$ is $\sum_{B=e^{\ell}, 0\leq \ell \leq [\log pq^{11/10}]}$.
We write
\bea \label{B91}
\mathcal {B}_1(p,q,t)=M_1+O(R_1),
\eea
where
\bea \label{M1}
M_1=\frac{qp^{\frac{1}{2}+it}}{12}\sum_{d\geq 1 \atop (d,q)=1}\frac{1}{d^{1+2it}}
\sum_{k =\diamond}\frac{\tau(k/p)}{k^{1+it}}W_t\left(\frac{4\pi^2 kd^2}{pq}\right)\sum_{d|(k/p,p)}\mu(d)=0,
\eea
since $(k/p,p)\neq 1$ if $k$ is a square, and
\bna
R_1=p^{\frac{1}{2}}\sum_{B}
\left(B^{\frac{7}{4}}q^{-\frac{1}{2}}+B^{\frac{39}{32}+\epsilon}q^{\frac{1}{4}}\right)
\sum_{d\geq 1 \atop (d,q)=1}\frac{1}{d}
\max_{B\leq k\leq eB}\frac{\tau(k/p)}{\sqrt{k}}\left|
W_t\left(\frac{4\pi^2 kd^2}{pq}\right)\right|.
\ena
To bound $R_1$, we note that there are $[\log pq^{11/10}]$ terms in the $B$-sum. Thus by Lemma 2.2,
\bea \label{R1}
R_1
&\ll_{\epsilon}&q^{\epsilon}\sqrt{p}\sum_{ B}\left(B^{\frac{5}{4}}q^{-\frac{1}{2}}+B^{\frac{23}{32}+\epsilon}q^{\frac{1}{4}}\right)
\exp\left(-C\sqrt{\frac{4\pi^2B}{pq}}\right)\nonumber\\
&\ll_{\epsilon}&q^{\epsilon}\left(p^{\frac{7}{4}}q^{\frac{3}{4}}+p^{\frac{39}{32}}q^{\frac{31}{32}}\right).
\eea
By \eqref{B91}-\eqref{R1}, for $p<q^{1/3}$, we have
\bea \label{B1pqt:result}
\mathcal {B}_1(p,q,t)\ll_{\epsilon} q^{\epsilon}\left(p^{\frac{7}{4}}q^{\frac{3}{4}}+p^{\frac{39}{32}}q^{\frac{31}{32}}\right)
\ll_{\epsilon} p^{\frac{39}{32}}q^{\frac{31}{32}+\epsilon}.
\eea

{\bf Evaluation of $\mathcal {B}_2(p,q,t)$}. Now we compute $\mathcal {B}_2(p,q,t)$ in (3.4).
We have
\bea \label{B2pqt}
\mathcal {B}_2(p,q,t)=\mathcal {B}_{21}(p,q,t)+\mathcal {B}_{22}(p,q,t),
\eea
where
\bea
\mathcal {B}_{21}(p,q,t)&=&\sum_{d\geq 1 \atop (d,q)=1}\frac{1}{d^{1+2it}}
\sum_{p|n}\frac{\tau(n)}{n^{\frac{1}{2}+it}}W_t\left(\frac{4\pi^2 nd^2}{q}\right)
\sum_{f\in \mathcal {S}_2^*(q)}\lambda_f(np),\\
\mathcal {B}_{22}(p,q,t)&=&\sum_{d\geq 1 \atop (d,q)=1}\frac{1}{d^{1+2it}}
\sum_{p|n}\frac{\tau(n)}{n^{\frac{1}{2}+it}}W_t\left(\frac{4\pi^2 nd^2}{q}\right)
\sum_{f\in \mathcal {S}_2^*(q)}\lambda_f\left(\frac{n}{p}\right).
\eea

We first compute $\mathcal {B}_{22}(p,q,t)$ in (3.11). We write
\bna
\mathcal {B}_{22}(p,q,t)=p^{-\frac{1}{2}-it}\sum_{d\geq 1 \atop (d,q)=1}\frac{1}{d^{1+2it}}
\sum_{k\geq 1}\frac{\tau(kp)}{k^{\frac{1}{2}+it}}W_t\left(\frac{4\pi^2 kpd^2}{q}\right)
Tr(T_k).
\ena
Again, by Lemma 2.2, we can assume that $kd^2\leq \exp([\log q^{11/10}])\leq q^{11/10}$.
Then we can apply Lemma 2.3 to obtain
\bea \label{B22pqt}
\mathcal {B}_{22}(p,q,t)=M_{22}+O(R_{22})
\eea
with
\bna
M_{22}=\frac{q}{12p^{\frac{1}{2}+it}}\sum_{d\geq 1 \atop (d,q)=1}\frac{1}{d^{1+2it}}
\sum_{k=\diamond}\frac{\tau(kp)}{k^{1+it}}W_t\left(\frac{4\pi^2 kpd^2}{q}\right)
\ena
and
\bna
R_{22}=\frac{1}{\sqrt{p}}\sum_B\left(B^{\frac{7}{4}}q^{-\frac{1}{2}}+B^{\frac{39}{32}+\epsilon}q^{\frac{1}{4}}\right)
\sum_{d\geq 1 \atop (d,q)=1}\frac{1}{d}\max_{B\leq k\leq eB}\frac{\tau(kp)}{\sqrt{k}}
\left|W_t\left(\frac{4\pi^2 kpd^2}{q}\right)\right|,
\ena
where the summation over $B$ is $\sum_{B=e^{\ell}, 0\leq \ell \leq [\log q^{11/10}]}$.
By Lemma 2.2 we have
\bea \label{R22}
R_{22}&\ll_{\epsilon}&\frac{q^{\epsilon}}{\sqrt{p}}
\sum_B\left(B^{\frac{5}{4}}q^{-\frac{1}{2}}+B^{\frac{23}{32}+\epsilon}q^{\frac{1}{4}}\right)
\exp\left(-C\sqrt{\frac{4\pi^2 Bp}{q}}\right)\nonumber\\
&\ll_{\epsilon}&q^{\epsilon}\left(p^{-\frac{7}{4}}q^{\frac{3}{4}}+p^{-\frac{39}{32}}q^{\frac{31}{32}}\right)
\ll_{\epsilon} q^{\frac{31}{32}+\epsilon}.
\eea

By (2.3), we have
\bea \label{M22}
M_{22}&=&\frac{q}{12p^{\frac{1}{2}+it}}\sum_{d\geq 1 \atop (d,q)=1}\frac{1}{d^{1+2it}}
\sum_{k=\diamond}\frac{\tau(kp)}{k^{1+it}}\frac{1}{2\pi i}\int\limits_{(2)}
\left(\frac{4\pi^2 kpd^2}{q}\right)^{-u}\Gamma(1+it+u)^2e^{u^2}\frac{\mbox{d}u}{u}\nonumber\\
&=&\frac{q}{12p^{\frac{1}{2}+it}}\frac{1}{2\pi i}\int\limits_{(2)}\zeta_q(1+2it+2u)
\left(\sum_{\ell\geq 1}\frac{\tau(p\ell^2)}{\ell^{2+2it+2u}}\right)
\left(\frac{4\pi^2 p}{q}\right)^{-u}\Gamma(1+it+u)^2e^{u^2}\frac{\mbox{d}u}{u},\nonumber\\
\eea
where $\zeta_q(s)=(1-q^{-s})\zeta(s)$.\\

Let $\ell=p^{\alpha}\ell'$ with $\alpha\geq 0$, $\ell'\geq 1$ and $(p,\ell')=1$. Then
$\tau(p\ell^2)=\tau(p^{2\alpha+1})\tau(\ell'^2)=(2\alpha+2)\tau(\ell'^2)$ and for $\Re(s)>1$,
\bna
\sum_{\ell\geq 1}\frac{\tau(p\ell^2)}{\ell^{s}}=
\sum_{\alpha\geq 0}\frac{2\alpha+2}{(p^{s})^\alpha}\sum_{(\ell',p)=1}\frac{\tau(\ell'^2)}{\ell'^{s}}=
\frac{2}{(1-p^{-s})^2}\frac{\zeta_p^3(s)}{\zeta_p(2s)}=\frac{2}{1+p^{-s}}\frac{\zeta^3(s)}{\zeta(2s)}.
\ena
Moving the line of integration in \eqref{M22} to $\Re(u)=-1/2+\epsilon$, we have
\bea \label{M222}
M_{22}=\frac{q}{12p^{\frac{1}{2}+it}}(M_{22}'+R_{22}'),
\eea
where
\bna
M_{22}'&=&\mbox{Res}_{u=0}\frac{2\zeta_q(1+2it+2u)}{1+p^{-2-2it-2u}}\frac{\zeta^3(2+2it+2u)}{\zeta(4+4it+4u)}
\left(\frac{4\pi^2 p}{q}\right)^{-u}\frac{\Gamma(1+it+u)^2}{u}e^{u^2},\\
R_{22}'&=&\frac{1}{2\pi i}\int\limits_{(-1/2+\epsilon)}
\frac{2\zeta_q(1+2it+2u)}{1+p^{-2-2it-2u}}\frac{\zeta^3(2+2it+2u)}{\zeta(4+4it+4u)}
\left(\frac{4\pi^2 p}{q}\right)^{-u}\Gamma(1+it+u)^2e^{u^2}\frac{\mbox{d}u}{u}.
\ena

By Stirling's formula and the convexity bound $\zeta(\sigma+i\tau)\ll (1+|\tau|)^{{(1-\sigma)/2}+\epsilon}$ for
$0<\sigma<1$, we have
\bea \label{R221}
R_{22}'&=&\frac{1}{2\pi}\int\limits_{-\infty}^{\infty}
\frac{2\zeta_q(2\epsilon+2it+2i\tau)}{1+p^{-1-2\epsilon-2it-2i\tau}}
\frac{\zeta^3(1+2\epsilon+2it+2i\tau)}{\zeta(2+4\epsilon+4it+4i\tau)}
\left(\frac{4\pi^2 p}{q}\right)^{\frac{1}{2}-\epsilon-i\tau}\nonumber\\
&&\times\Gamma\left(\frac{1}{2}+\epsilon+it+i\tau\right)^2
e^{(-\frac{1}{2}+\epsilon+i\tau)^2}
\frac{\mbox{d}\tau}{-\frac{1}{2}+\epsilon+i\tau}\nonumber\\
&=&\frac{1}{2\pi}\int\limits_{-\infty}^{\infty}
\frac{2\zeta_q(2\epsilon+2i\tau)}{1+p^{-1-2\epsilon-2i\tau}}
\frac{\zeta^3(1+2\epsilon+2i\tau)}{\zeta(2+4\epsilon+4i\tau)}
\left(\frac{4\pi^2 p}{q}\right)^{\frac{1}{2}-\epsilon-i(\tau-t)}\nonumber\\
&&\times\Gamma\left(\frac{1}{2}+\epsilon+i\tau\right)^2
e^{(-\frac{1}{2}+\epsilon+i(\tau-t))^2}
\frac{\mbox{d}\tau}{-\frac{1}{2}+\epsilon+i(\tau-t)}\nonumber\\
&\ll&\left(\frac{p}{q}\right)^{\frac{1}{2}-\epsilon}\int\limits_{-\infty}^{\infty}
(1+|\tau|)^{{\frac{1}{2}(1-2\epsilon)}+\epsilon}
\left|\Gamma\left(\frac{1}{2}+\epsilon+i\tau\right)\right|^2\mbox{d}\tau\nonumber\\
&\ll&\left(\frac{p}{q}\right)^{\frac{1}{2}-\epsilon}\left(1+\int\limits_{|\tau|>1}
|\tau|^{{\frac{1}{2}(1-2\epsilon)}+\epsilon}|\tau|^{2\epsilon}e^{-\pi |\tau|}\mbox{d}\tau\right)\nonumber\\
&\ll&\left(\frac{p}{q}\right)^{\frac{1}{2}-\epsilon}.
\eea

We compute $M_{22}'$ according as $t=0$ and $t\neq 0$.

\noindent Case I: If $t=0$, there is a double pole at $u=0$. We have
\bna
M_{22}'&=&\lim_{u\rightarrow 0} \frac{\mbox{d}}{\mbox{d}u}
\left\{\frac{2u\zeta_q(1+2u)}{1+p^{-2-2u}}\frac{\zeta^3(2+2u)}{\zeta(4+4u)}
\left(\frac{4\pi^2 p}{q}\right)^{-u}\Gamma(1+u)^2e^{u^2}\right\}.\nonumber
\ena
Since $\zeta_q(s)=(1-q^{-s})\zeta(s)$, we have
\bna
M_{22}'=\lim_{u\rightarrow 0} \frac{\mbox{d}}{\mbox{d}u}
\left\{g_1(u)g_2(u)\right\}=g_1'(0)g_2(0)+g_1(0)g_2'(0),
\ena
where
\bna
&&g_1(u)=u\zeta\left({1+2u}\right),\\
&&g_2(u)=2\frac{1-q^{-1-2u}}{1+p^{-2-2u}}\frac{\zeta^3\left(2+2u\right)}{\zeta\left(4+4u\right)}
\left(\frac{4\pi^2 p}{q}\right)^{-u}\Gamma(1+u)^2e^{u^2}.
\ena
By the Laurent expansion of $\zeta(s)$,
\bna
\zeta(1+2u)=\frac{1}{2u}+\gamma+\sum_{n=1}^{\infty}(-1)^n\frac{\gamma_n}{n!}(2u)^n,
\ena
where $\gamma$ is Euler's constant and $\gamma_n$ are stieltjes constants.
Thus $g_1(0)=\frac{1}{2}$ and
\bna
g'_1(u)=\gamma+\sum_{n=1}^{\infty}\frac{(-2)^n(n+1)\gamma_n}{n!}u^n,
\ena
which in turn implies $g_1'(0)=\gamma$. For $g_2(u)$, we note that
\bna
g_2(0)&=&2\frac{1-q^{-1}}{1+p^{-2}}\frac{\zeta^3\left(2\right)}{\zeta\left(4\right)},\\
g'_2(0)&=&\left(\frac{2q^{-1}\log{q}}{1-q^{-1}}
+\frac{2p^{-2}\log{p}}{1+p^{-2}}
+\frac{6\zeta'(2)}{\zeta(2)}
-\frac{4\zeta'(4)}{\zeta{(4)}}
-\log{\frac{4\pi^2p}{q}}
+2\Gamma'(1)\right)g_2(0)\\
&=&\left(\frac{6\zeta'(2)}{\zeta(2)}
-\frac{4\zeta'(4)}{\zeta{(4)}}
+\frac{2\log{p}}{1+p^2}
-\log{\frac{4\pi^2p}{q}}
-2\gamma\right)g_2(0)
+O\left(q^{-1}\log q\right),
\ena
where we have used the fact that $\Gamma'(1)=-\gamma$.

Assembling the above computations, we get
\bna
M_{22}'&=&\frac{1-q^{-1}}{1+p^{-2}}\frac{\zeta^3(2)}{\zeta(4)}
\left(\log{\frac{q}{4\pi^2p}}
+\frac{2\log{p}}{1+p^2}
+6\frac{\zeta'(2)}{\zeta(2)}
-4\frac{\zeta'(4)}{\zeta(4)}\right)+O\left(q^{-1}\log q\right).
\ena
This result combined with \eqref{M222} and \eqref{R221} yields
\bea \label{M2222}
M_{22}=\frac{1-q^{-1}}{1+p^{-2}}
\frac{\zeta^3(2)}{12\zeta(4)}\frac{q}{\sqrt{p}}
\left(\log\frac{q}{4\pi^2p}
+\frac{2\log{p}}{1+p^2}
+6\frac{\zeta'(2)}{\zeta(2)}
-4\frac{\zeta'(4)}{\zeta(4)}\right)
+O_\epsilon\left(q^{\frac{1}{2}+\epsilon}\right).
\eea
By \eqref{B22pqt}, \eqref{R22} and \eqref{M2222}, we conclude that
\bea \label{B22pq0}
\mathcal {B}_{22}(p,q,0)=\frac{1-q^{-1}}{1+p^{-2}}
\frac{\zeta^3(2)}{12\zeta(4)}\frac{q}{\sqrt{p}}
\left(\log\frac{q}{4\pi^2p}
+\frac{2\log{p}}{1+p^2}
+6\frac{\zeta'(2)}{\zeta(2)}
-4\frac{\zeta'(4)}{\zeta(4)}\right)
+O_{\epsilon}\left(q^{\frac{31}{32}+\epsilon}\right).
\eea

\medskip

\noindent Case II: If $t\neq 0$, there are simple poles at $u=0$ and $u=-it$. Then we have
\bna
M_{22}'&=&\lim_{u\rightarrow -it} \frac{2(u+it)\zeta_q(1+2it+2u)}{1+p^{-2-2it-2u}}\frac{\zeta^3(2+2it+2u)}{\zeta(4+4it+4u)}
\left(\frac{4\pi^2 p}{q}\right)^{-u}\frac{\Gamma(1+it+u)^2}{u}e^{u^2}\nonumber\\
&&+\lim_{u\rightarrow 0} \frac{2\zeta_q(1+2it+2u)}{1+p^{-2-2it-2u}}\frac{\zeta^3(2+2it+2u)}{\zeta(4+4it+4u)}
\left(\frac{4\pi^2 p}{q}\right)^{-u}\Gamma(1+it+u)^2e^{u^2}\nonumber\\
&=&\frac{1-q^{-1}}{1+p^{-2}}\frac{\zeta^3(2)}{\zeta(4)}\left(\frac{4\pi^2p}{q}\right)^{it}\frac{-e^{-t^2}}{it}
+\frac{2\zeta_q(1+2it)}{1+p^{-2-2it}}\frac{\zeta^3(2+2it)}{\zeta(4+4it)}\Gamma(1+it)^2.
\ena
Then by \eqref{M222} and \eqref{R221}, we obtain
\bna
M_{22}=\frac{q}{12p^{\frac{1}{2}+it}}\left(\frac{1-q^{-1}}{1+p^{-2}}\frac{\zeta^3(2)}{\zeta(4)}\left(\frac{4\pi^2p}{q}\right)^{it}\frac{-e^{-t^2}}{it}
+\frac{2\zeta_q(1+2it)}{1+p^{-2-2it}}\frac{\zeta^3(2+2it)}{\zeta(4+4it)}\Gamma(1+it)^2\right)+
O\left(q^{\frac{1}{2}+\epsilon}\right).
\ena
Moreover, by \eqref{B22pqt} and \eqref{R22}, we have
\bea \label{B22pqt:result}
\mathcal {B}_{22}(p,q,t)&=&\frac{q}{12p^{\frac{1}{2}+it}}
\left(\frac{1-q^{-1}}{1+p^{-2}}\frac{\zeta^3(2)}{\zeta(4)}\left(\frac{4\pi^2p}{q}\right)^{it}\frac{-e^{-t^2}}{it}
+\frac{2\zeta_q(1+2it)}{1+p^{-2-2it}}\frac{\zeta^3(2+2it)}{\zeta(4+4it)}\Gamma(1+it)^2\right)\nonumber\\
&&+O\left(q^{\frac{31}{32}+\epsilon}\right).
\eea

Next we compute $\mathcal {B}_{21}(p,q,t)$ in (3.10).
We write
\bna
\mathcal {B}_{21}(p,q,t)&=&p^{\frac{1}{2}+it}\sum_{d\geq 1 \atop (d,q)=1}\frac{1}{d^{1+2it}}
\sum_{p|\frac{m}{p}}\frac{\tau(m/p)}{m^{\frac{1}{2}+it}}W_t\left(\frac{4\pi^2 md^2}{pq}\right)
\sum_{f\in \mathcal {S}_2^*(q)}\lambda_f(m)\\
&=&p^{\frac{1}{2}+it}\sum_{d\geq 1 \atop (d,q)=1}\frac{1}{d^{1+2it}}
\sum_{m\geq p^2\atop p^2|m}\frac{\tau(m/p)}{m^{\frac{1}{2}+it}}W_t\left(\frac{4\pi^2 md^2}{pq}\right)T_r(T_m).
\ena
Again, by Lemma 2.2, we can assume that $md^2\leq \exp([\log pq^{11/10}])$.
Then we can apply Lemma 2.3 to obtain
\bea \label{B21pqt}
\mathcal {B}_{21}(p,q,t)=M_{21}+O(R_{21})
\eea
with
\bna
M_{21}=\frac{qp^{\frac{1}{2}+it}}{12}\sum_{d\geq 1 \atop (d,q)=1}\frac{1}{d^{1+2it}}\sum_{m=\diamond \atop p^2|m}\frac{\tau(m/p)}{m^{1+it}}W_t\left(\frac{4\pi^2 md^2}{pq}\right)
\ena
and
\bna
R_{21}=p^{\frac{1}{2}}\sum_{B}
\left(B^{\frac{7}{4}}q^{-\frac{1}{2}}+B^{\frac{39}{32}+\epsilon}q^{\frac{1}{4}}\right)
\sum_{d\geq 1 \atop (d,q)=1}\frac{1}{d}
\max_{B\leq m\leq eB \atop p^2\mid m}\frac{\tau(m/p)}{\sqrt{m}}\left|
W_t\left(\frac{4\pi^2 md^2}{pq}\right)\right|,
\ena
where the summation over $B$ is $\sum_{B=e^{\ell}, [\log p^2]\leq \ell \leq [\log pq^{11/10}]}$.

Similar as the proof of (3.8), we have for $p<q^{1/3}$,
\bea \label{R21}
R_{21}
&\ll_{\epsilon}&p^{\frac{39}{32}}q^{\frac{31}{32}+\epsilon}.
\eea
For $M_{21}$, by (2.3) and \eqref{M22},
\bea \label{M21}
M_{21}&=&\frac{qp^{\frac{1}{2}+it}}{12}\sum_{d\geq 1 \atop (d,q)=1}\frac{1}{d^{1+2it}}
\sum_{m=\diamond \atop p^2 \mid m}\frac{\tau(m/p)}{m^{1+it}}\frac{1}{2\pi i}\int\limits_{(2)}
\left(\frac{4\pi^2 md^2}{pq}\right)^{-u}\Gamma(1+it+u)^2e^{u^2}\frac{\mbox{d}u}{u}\nonumber\\
&=&\frac{qp^{\frac{1}{2}+it}}{12}\sum_{d\geq 1 \atop (d,q)=1}\frac{1}{d^{1+2it}}
\sum_{\ell \geq 1}\frac{\tau(p\ell^2)}{(p\ell)^{2+2it}}\frac{1}{2\pi i}\int\limits_{(2)}
\left(\frac{4\pi^2 (p\ell)^2d^2}{pq}\right)^{-u}\Gamma(1+it+u)^2e^{u^2}\frac{\mbox{d}u}{u}\nonumber\\
&=&\frac{q}{12p^{\frac{3}{2}+it}}\frac{1}{2\pi i}\int\limits_{(2)}\zeta_q(1+2it+2u)
\left(\sum_{\ell\geq 1}\frac{\tau(p\ell^2)}{\ell^{2+2it+2u}}\right)
\left(\frac{4\pi^2 p}{q}\right)^{-u}\Gamma(1+it+u)^2e^{u^2}\frac{\mbox{d}u}{u}\nonumber\\
&=&\frac{1}{p}M_{22}.
\eea
By \eqref{B21pqt}-\eqref{M21},
\bea \label{B21pqt1}
\mathcal {B}_{21}(p,q,t)&=&\frac{1}{p}M_{22}+O_\epsilon\left(p^{\frac{39}{32}}q^{\frac{31}{32}+\epsilon}\right).
\eea

Consequently, if $ t=0 $ ,by \eqref{B22pq0} and \eqref{B21pqt1},
\bea \label{B2pq0:result}
\mathcal {B}_2(p,q,0)&=&\mathcal{B}_{21}(p,q,0)+\mathcal {B}_{22}(p,q,0)\nonumber\\
&=&\left(1+\frac{1}{p}\right)\frac{1-q^{-1}}{1+p^{-2}}
\frac{\zeta^3(2)}{12\zeta(4)}\frac{q}{\sqrt{p}}
\left(\log\frac{q}{4\pi^2p}
+\frac{2\log{p}}{1+p^2}
+6\frac{\zeta'(2)}{\zeta(2)}
-4\frac{\zeta'(4)}{\zeta(4)}\right)\nonumber\\
&&+O_\epsilon\left(p^{\frac{39}{32}}q^{\frac{31}{32}+\epsilon}\right).
\eea

If $t\neq 0$ , then by \eqref{B22pqt:result} and \eqref{B21pqt1}, we can obtain
\bea \label{B2pqt:result}
\mathcal {B}_{2}(p,q,t)&=&\mathcal{B}_{21}(p,q,t)+\mathcal{B}_{22}(p,q,t)\nonumber\\
&=&\frac{q(p+1)}{12p^{\frac{3}{2}+it}}
\left(\frac{1-q^{-1}}{1+p^{-2}}\frac{\zeta^3(2)}{\zeta(4)}\left(\frac{4\pi^2p}{q}\right)^{it}\frac{-e^{-t^2}}{it}
+\frac{2\zeta_q(1+2it)}{1+p^{-2-2it}}\frac{\zeta^3(2+2it)}{\zeta(4+4it)}\Gamma(1+it)^2\right)\nonumber\\
&&+O_\epsilon\left(p^{\frac{39}{32}}q^{\frac{31}{32}+\epsilon}\right).
\eea
{\bf Conclusion}. By (3.2), \eqref{B1pqt:result} and \eqref{B2pq0:result}, we have for $t=0$ ,
\bna
\mathcal {B}(p,q,0)=\left(1+\frac{1}{p}\right)
\frac{1-q^{-1}}{1+p^{-2}}
\frac{\zeta^3(2)}{12\zeta(4)}\frac{q}{\sqrt{p}}
\left(\log\frac{q}{4\pi^2p}
+\frac{2\log{p}}{1+p^2}
+6\frac{\zeta'(2)}{\zeta(2)}
-4\frac{\zeta'(4)}{\zeta(4)}\right)+O_{\epsilon}\left(p^{\frac{39}{32}}q^{\frac{31}{32}+\epsilon}\right),
\ena
and by (3.1),
\bna
\mathcal {A}(p,q,0)&=&\left(1+\frac{1}{p}\right)\frac{1-q^{-1}}{1+p^{-2}}
\frac{\zeta^3(2)}{6\zeta(4)}\frac{q}{\sqrt{p}}
\left(\log\frac{q}{4\pi^2p}
+\frac{2\log{p}}{1+p^2}
+6\frac{\zeta'(2)}{\zeta(2)}-4\frac{\zeta'(4)}{\zeta(4)}\right)
+O_{\epsilon}\left(p^{\frac{39}{32}}q^{\frac{31}{32}+\epsilon}\right),
\ena
which is an asymptotic formula for $p<q^{\frac{1}{55}-\epsilon}$.

Let $0\neq t\in \mathbb{R}$ be fixed.
By (3.2), \eqref{B1pqt:result} and \eqref{B2pqt:result}, we have
\bea \label{Bpqt:result}
\mathcal {B}(p,q,t)&=&\frac{q(p+1)}{12p^{\frac{3}{2}+it}}
\left(\frac{1-q^{-1}}{1+p^{-2}}\frac{\zeta^3(2)}{\zeta(4)}\left(\frac{4\pi^2p}{q}\right)^{it}\frac{-e^{-t^2}}{it}
+\frac{2\zeta_q(1+2it)}{1+p^{-2-2it}}\frac{\zeta^3(2+2it)}{\zeta(4+4it)}\Gamma(1+it)^2\right)\nonumber\\
&&+O_\epsilon\left(p^{\frac{39}{32}}q^{\frac{31}{32}+\epsilon}\right).
\eea
By (3.1) and \eqref{Bpqt:result}, we conclude that
\bna
\mathcal {A}(p,q,t)&=&
\frac{\Gamma(1-it)^2}{\Gamma(1+it)^2}\left(\frac{q}{4\pi^2}\right)^{-2it}
\frac{\zeta_q(1-2it)}{6(1+p^{-2+2it})}\frac{\zeta^3(2-2it)}{\zeta(4-4it)}\frac{q(p+1)}{p^{\frac{3}{2}-it}}\nonumber\\
&&+\frac{\zeta_q(1+2it)}{6(1+p^{-2-2it})}\frac{\zeta^3(2+2it)}{\zeta(4+4it)}\frac{q(p+1)}{p^{\frac{3}{2}+it}}
+O_{t,\epsilon}\left(p^{\frac{39}{32}}q^{\frac{31}{32}+\epsilon}\right),
\ena
which is an asymptotic formula for $p<q^{\frac{1}{55}-\epsilon}$.

\noindent {\bf 4 \quad Proof of Theorem 1.2}
\setcounter{section}{4} \setcounter{equation}{0}

\medskip
In this section, we evaluate the sum
\bna
\mathcal {A}(p^2,q,t)=\sum_{f\in S_2^*(q)} L\left(\frac{1}{2}+it,f\right)^2\lambda_f(p^2).
\ena

The proof is similar as that of Thoerem 1.1, so we will be brief.
By the approximate functional equation in (2.2), we have
\bea \label{Ap2qt}
\mathcal {A}(p^2,q,t)
=\Gamma(1+it)^{-2}\mathcal {B}(p^2,q,t)+\Gamma(1+it)^{-2}\left(\frac{q}{4\pi^2}\right)^{-2it}\overline{\mathcal {B}(p^2,q,t)},
\eea
where
\bna
\mathcal {B}(p^2,q,t)=\sum_{d\geq 1 \atop (d,q)=1}\frac{1}{d^{1+2it}}
\sum_{n\geq 1}\frac{\tau(n)}{n^{\frac{1}{2}+it}}W_t\left(\frac{4\pi^2 nd^2}{q}\right)
\sum_{f\in \mathcal {S}_2^*(q)}\lambda_f(n)\lambda_f(p^2).
\ena
By (2.1),
\bna
\lambda_f(n)\lambda_f(p^2)=\sum_{d|(n,p^2) \atop (d,q)=1}\lambda_f\left(\frac{np^2}{d^2}\right)
=\left\{\begin{array}{lll}
 \lambda_f(np^2), & \mbox {if $p \nmid n$},\\
 \lambda_f(np^2)+\lambda_f(n), & \mbox {if $p \mid n, p^2 \nmid n$},\\
\lambda_f(np^2)+\lambda_f(n)+\lambda_f\left(\frac{n}{p^2}\right), &\mbox {if $p \mid n, p^2 \mid n$.}\end{array}\right.
\ena
Correspondingly,
\bea \label{BB}
\mathcal {B}(p^2,q,t)=\mathcal {B}_1(p^2,q,t)
+\mathcal {B}_2(p^2,q,t)
+\mathcal {B}_3(p^2,q,t),
\eea
where
\bna
\mathcal {B}_1(p^2,q,t)&=&\sum_{d\geq 1 \atop (d,q)=1}\frac{1}{d^{1+2it}}
\sum_{n\geq1}\frac{\tau(n)}{n^{\frac{1}{2}+it}}W_t\left(\frac{4\pi^2 nd^2}{q}\right)
\sum_{f\in \mathcal {S}_2^*(q)}\lambda_f\left(np^2\right),\\
\mathcal {B}_2(p^2,q,t)&=&\sum_{d\geq 1 \atop (d,q)=1}\frac{1}{d^{1+2it}}
\sum_{p|n}\frac{\tau(n)}{n^{\frac{1}{2}+it}}W_t\left(\frac{4\pi^2 nd^2}{q}\right)
\sum_{f\in \mathcal {S}_2^*(q)}\lambda_f(n),\\
\mathcal {B}_3(p^2,q,t)&=&\sum_{d\geq 1 \atop (d,q)=1}\frac{1}{d^{1+2it}}
\sum_{p^2|n}\frac{\tau(n)}{n^{\frac{1}{2}+it}}W_t\left(\frac{4\pi^2 nd^2}{q}\right)
\sum_{f\in \mathcal {S}_2^*(q)}\lambda_f\left(\frac{n}{p^2}\right).
\ena

\medskip
{\bf Evaluation of $\mathcal {B}_1(p^2,q,t)$}. We have
\bea \label{B1}
\mathcal {B}_1(p^2,q,t)
&=&p^{1+2it}\sum_{d\geq 1 \atop (d,q)=1}\frac{1}{d^{1+2it}}
\sum_{k\geq1 \atop p^2\mid k}\frac{\tau(k/p^2)}{k^{\frac{1}{2}+it}}
W_t\left(\frac{4\pi^2 kd^2}{p^2q}\right)Tr\left(T_k\right).
\eea
By Lemma 2.2, we can assume that $kd^2\leq \exp([\log p^2q^{11/10}])$.
Then by Lemma 2.3,
\bea \label{B11}
\mathcal {B}_1(p^2,q,t)&=&\frac{qp^{1+2it}}{12}\sum_{d\geq 1 \atop (d,q)=1}\frac{1}{d^{1+2it}}
\sum_{k =\diamond \atop p^2\mid k}\frac{\tau(k/p^2)}{k^{1+it}}W_t\left(\frac{4\pi^2 kd^2}{p^2q}\right)\nonumber\\
&+&O\left(p\sum_{B}
\left(B^{\frac{7}{4}}q^{-\frac{1}{2}}+B^{\frac{39}{32}+\epsilon}q^{\frac{1}{4}}\right)
\sum_{d\geq 1 \atop (d,q)=1}\frac{1}{d}
\max_{B\leq k\leq eB \atop p^2 \mid k}\frac{\tau(k/p^2)}{\sqrt{k}}\left|
W_t\left(\frac{4\pi^2 kd^2}{p^2q}\right)\right|\right)\nonumber\\
&:=&\Delta_1+E_1,
\eea
say, where the summation over $B$ is $\sum_{B=e^{\ell}, 0\leq \ell \leq [\log p^2q^{11/10}]}$.
To bound $E_1$, by Lemma 2.2,
\bna
E_1
&\ll_{\epsilon}&pq^{\epsilon}\sum_{ B}\left(B^{\frac{5}{4}}q^{-\frac{1}{2}}+B^{\frac{23}{32}+\epsilon}q^{\frac{1}{4}}\right)
\exp\left(-C\sqrt{\frac{4\pi^2 B}{p^2q}}\right)\nonumber\\
&\ll_{\epsilon}&pq^\epsilon\left(
p^{\frac{5}{2}}q^{\frac{3}{4}+\epsilon}+p^{\frac{23}{16}+\epsilon}q^{\frac{31}{32}+\epsilon}\right)\nonumber\\
&\ll_{\epsilon}&p^{\frac{7}{2}}q^{\frac{3}{4}+\epsilon}+p^{\frac{39}{16}+\epsilon}q^{\frac{31}{32}+\epsilon}.
\ena
Then by assuming $p<q^{1/6}$, we get
\bea \label{E1}
E_1\ll_{\epsilon}
p^{\frac{39}{16}+\epsilon}q^{\frac{31}{32}+\epsilon}.
\eea
By (2.3), we have
\bea \label{D1}
\Delta_1&=&\frac{qp^{1+2it}}{12}\sum_{d\geq 1 \atop (d,q)=1}\frac{1}{d^{1+2it}}
\sum_{k=\diamond \atop p^2\mid k}\frac{\tau (k/p^2)}{k^{1+it}}\frac{1}{2\pi i}\int\limits_{(2)}
\left(\frac{4\pi^2 kd^2}{p^2q}\right)^{-u}\Gamma(1+it+u)^2e^{u^2}
\frac{\mbox{d}u}{u}\nonumber\\
&=&\frac{q}{12p}\frac{1}{2\pi i}\int\limits_{(2)}\zeta_q(1+2it+2u)
\left(\sum_{\ell \geq1}\frac{\tau(\ell^2)}{\ell^{2+2it+2u}}\right)
\left(\frac{4\pi^2}{q}\right)^{-u}
\Gamma(1+it+u)^2e^{u^2}\frac{\mbox{d}u}{u}\nonumber\\
&=&\frac{q}{12p}\frac{1}{2\pi i}\int\limits_{(2)}
\frac{\zeta_q(1+2it+2u)\zeta^3(2+2it+2u)}{\zeta(4+4it+4u)}
\left(\frac{4\pi^2}{q}\right)^{-u}
\Gamma(1+it+u)^2e^{u^2}\frac{\mbox{d}u}{u}.
\eea
Similar as the evaluation of $M_{22}$, we move the line of integration in \eqref{D1} to $\Re(u)=-1/2+\epsilon$ to get
\bea \label{D11}
\Delta_1=\frac{q}{12p}(\Delta_1'+E_1'),
\eea
where
\bna
\Delta_1'&=&\mbox{Res}_{u=0}\frac{
\zeta_q(1+2it+2u)\zeta^3(2+2it+2u)}
{\zeta(4+4it+4u)}
\left(\frac{4\pi^2}{q}\right)^{-u}
\frac{\Gamma(1+it+u)^2}{u}e^{u^2},\\
E_1'&=&\frac{1}{2\pi i}\int\limits_{(-1/2+\epsilon)}
\frac{\zeta_q(1+2it+2u)\zeta^3(2+2it+2u)}
{\zeta(4+4it+4u)}
\left(\frac{4\pi^2}{q}\right)^{-u}
\Gamma(1+it+u)^2e^{u^2}\frac{\mbox{d}u}{u}.
\ena
By Stirling's formula and the convexity bound $\zeta(\sigma+i\tau)\ll (1+|\tau|)^{{(1-\sigma)/2}+\epsilon}$ for $0<\sigma<1$,
\bea \label{E11}
E_1'
&\ll&q^{-\frac{1}{2}+\epsilon}
\int\limits_{-\infty}^{\infty}
(1+|\tau|)^{{\frac{1}{2}(1-2\epsilon)}+\epsilon}
\left|\Gamma\left(\frac{1}{2}+\epsilon+i\tau\right)\right|^2\mbox{d}\tau\nonumber\\
&\ll&q^{-\frac{1}{2}+\epsilon}.
\eea

If $t=0$, there is a double pole at $u=0$ and
\bea \label{D111}
\Delta_1'&=&\lim_{u\rightarrow 0} \frac{\mbox{d}}{\mbox{d}u}
\left\{\frac{u\zeta_q(1+2u)\zeta^3(2+2u)}{\zeta(4+4u)}
\left(\frac{4\pi^2}{q}\right)^{-u}
\Gamma(1+u)^2e^{u^2}\right\}\nonumber\\
&=&(1-q^{-1})
\frac{\zeta^3(2)}{2\zeta(4)}\left(\log{\frac{q}{4\pi^2}}
+6\frac{\zeta'(2)}{\zeta(2)}-4\frac{\zeta'(4)}{\zeta(4)}\right)
+O\left(q^{-1}\log q\right).
\eea
Then by \eqref{D11}-\eqref{D111}, we obtain
\bna
  \Delta_1=\frac{q-1}{24p}
\frac{\zeta^3(2)}{\zeta(4)}
\left(\log\frac{q}{4\pi^2}
+\frac{6\zeta'(2)}{\zeta(2)}-\frac{4\zeta'(4)}{\zeta(4)}\right)
+O\left(p^{-1}q^{\frac{1}{2}+\epsilon}\right).
\ena
Moreover, by \eqref{B11} and \eqref{E1}, we have
\bea \label{B10}
\mathcal {B}_1(p^2,q,0)
=\frac{q-1}{24p}
\frac{\zeta^3(2)}{\zeta(4)}
\left(\log\frac{q}{4\pi^2}
+\frac{6\zeta'(2)}{\zeta(2)}-\frac{4\zeta'(4)}{\zeta(4)}\right)
+O\left(p^{\frac{39}{16}+\epsilon}q^{\frac{31}{32}+\epsilon}\right).
\eea

If $t\neq 0$, there are simple poles at $u=0$ and $u=-it$ and
\bea \label{D5}
\Delta_1'&=&\lim_{u\rightarrow -it}
\frac{(u+it)\zeta_q(1+2it+2u)\zeta^3(2+2it+2u)}{\zeta(4+4it+4u)}
\left(\frac{4\pi^2}{q}\right)^{-u}\frac{\Gamma(1+it+u)^2}{u}e^{u^2}\nonumber\\
&+&\lim_{u\rightarrow 0}
\frac{\zeta_q(1+2it+2u)\zeta^3(2+2it+2u)}{\zeta(4+4it+4u)}
\left(\frac{4\pi^2}{q}\right)^{-u}\Gamma(1+it+u)^2e^{u^2}\nonumber\\
&=&\frac{(1-q^{-1})\zeta^3(2)}{2\zeta(4)}\left(\frac{4\pi^2}{q}\right)^{it}\frac{-e^{-t^2}}{it}
+\frac{\zeta_q(1+2it)\zeta^3(2+2it)}{\zeta(4+4it)}\Gamma(1+it)^2.
\eea
By \eqref{D11}, \eqref{E11} and \eqref{D5}, we obtain,
\bna
\Delta_1&=&
\frac{q}{12p}\left(
\frac{(1-q^{-1})\zeta^3(2)}{2\zeta(4)}\left(\frac{4\pi^2}{q}\right)^{it}\frac{-e^{-t^2}}{it}
+\frac{\zeta_q(1+2it)\zeta^3(2+2it)}{\zeta(4+4it)}\Gamma(1+it)^2\right)
+O\left(p^{-1}q^{\frac{1}{2}+\epsilon}\right).
\ena
Moreover, by \eqref{B11} and \eqref{E1}, we have
\bea \label{B1p2qt}
\mathcal {B}_1(p^2,q,t)&=&
\frac{q}{12p}\left(
\frac{(1-q^{-1})\zeta^3(2)}{2\zeta(4)}\left(\frac{4\pi^2}{q}\right)^{it}\frac{-e^{-t^2}}{it}
+\frac{\zeta_q(1+2it)\zeta^3(2+2it)}{\zeta(4+4it)}\Gamma(1+it)^2\right)\nonumber\\
&+&O\left(p^{\frac{39}{16}+\epsilon}q^{\frac{31}{32}+\epsilon}\right).
\eea

\medskip

{\bf Evaluation of $\mathcal {B}_3(p^2,q,t)$}.
We write
\bna
\mathcal {B}_3(p^2,q,t)
&=&\frac{1}{p^{1+2it}}\sum_{d\geq 1 \atop (d,q)=1}\frac{1}{d^{1+2it}}
\sum_{m\geq 1}\frac{\tau(mp^2)}{{m}^{\frac{1}{2}+it}}W_t\left(\frac{4\pi^2 mp^2d^2}{q}\right)Tr(T_m).
\ena
By Lemma 2.2 and 2.3,
\bea \label{B3}
\mathcal {B}_3(p^2,q,t)&=&
\frac{q}{12p^{1+2it}}\sum_{d\geq 1 \atop (d,q)=1}\frac{1}{d^{1+2it}}
\sum_{m=\diamond}\frac{\tau(mp^2)}{m^{1+it}}W_t\left(\frac{4\pi^2 mp^2d^2}{q}\right)\nonumber\\
&+&O\left(\frac{1}{p}\sum_B\left(B^{\frac{7}{4}}q^{-\frac{1}{2}}+B^{\frac{39}{32}+\epsilon}q^{\frac{1}{4}}\right)
\sum_{d\geq 1 \atop (d,q)=1}\frac{1}{d}\max_{B\leq m\leq eB}\frac{\tau(mp^2)}{\sqrt{m}}
\left|W_t\left(\frac{4\pi^2 mp^2d^2}{q}\right)\right|\right)\nonumber\\
&:=&\Delta_3+E_3,
\eea
say, where the summation over $B$ is $\sum_{B=e^{\ell}, 0\leq \ell \leq [\log q^{11/10}]}$.
By Lemma 2.2 we have
\bea \label{E3}
E_3&\ll_{\epsilon}&\frac{1}{p}q^{\epsilon}
\sum_B\left(B^{\frac{5}{4}}q^{-\frac{1}{2}}+B^{\frac{23}{32}+\epsilon}q^{\frac{1}{4}}\right)
\exp\left(-C\sqrt{\frac{4\pi^2 Bp^2}{q}}\right)\nonumber\\
&\ll_{\epsilon}&q^{\frac{31}{32}+\epsilon}.
\eea

By (2.3), we have
\bea \label{D31}
\Delta_3&=&\frac{q}{12p^{1+2it}}\sum_{d\geq 1 \atop (d,q)=1}\frac{1}{d^{1+2it}}
\sum_{\ell \geq 1}\frac{\tau( p^2\ell^2)}{\ell^{2+2it}}
\frac{1}{2\pi i}\int\limits_{(2)}
\left(\frac{4\pi^2 \ell^2p^2d^2}{q}\right)^{-u}\Gamma(1+it+u)^2e^{u^2}\frac{\mbox{d}u}{u}\nonumber\\
&=&\frac{q}{12p^{1+2it}}\frac{1}{2\pi i}\int\limits_{(2)}\zeta_q(1+2it+2u)
\left(\sum_{\ell\geq 1}\frac{\tau(p^2\ell^2)}{\ell^{2+2it+2u}}\right)
\left(\frac{4\pi^2 p^2}{q}\right)^{-u}\Gamma(1+it+u)^2e^{u^2}\frac{\mbox{d}u}{u}\nonumber\\
&=&\frac{q}{12p^{1+2it}}\frac{1}{2\pi i}\int\limits_{(2)}
\frac{(3-p^{-2-2it-2u})\zeta_q(1+2it+2u)\zeta^3(2+2it+2u)}{(1+p^{-2-2it-2u})\zeta(4+4it+4u)}
\left(\frac{4\pi^2p^2}{q}\right)^{-u}
\Gamma(1+it+u)^2e^{u^2}\frac{\mbox{d}u}{u}.\nonumber\\
\eea
Moving the line of integration in \eqref{D31} to $\Re(u)=-1/2+\epsilon$, we have
\bea \label{D32}
\Delta_3=\frac{q}{12p^{1+2it}}(\Delta_3'+E_3'),
\eea
where
\bna
\Delta_3'&=&\mbox{Res}_{u=0}\frac{(3-p^{-2-2it-2u})\zeta_q(1+2it+2u)\zeta^3(2+2it+2u)}
{(1+p^{-2-2it-2u})\zeta(4+4it+4u)}
\left(\frac{4\pi^2 p^2}{q}\right)^{-u}\frac{\Gamma(1+it+u)^2}{u}e^{u^2},\\
E_3'&=&\frac{1}{2\pi i}\int\limits_{(-1/2+\epsilon)}
\frac{(3-p^{-2-2it-2u})\zeta_q(1+2it+2u)\zeta^3(2+2it+2u)}
{(1+p^{-2-2it-2u})\zeta(4+4it+4u)}
\left(\frac{4\pi^2 p^2}{q}\right)^{-u}\Gamma(1+it+u)^2e^{u^2}\frac{\mbox{d}u}{u}.
\ena

By Stirling's formula and the convexity bound $\zeta(\sigma+i\tau)\ll (1+|\tau|)^{{(1-\sigma)/2}+\epsilon}$ for
$0<\sigma<1$,
\bea \label{E31}
E_3'
&\ll&\left(\frac{p^2}{q}\right)^{\frac{1}{2}-\epsilon}.
\eea

If $t=0$, there is a double pole at $u=0$, then
\bea \label{D33}
\Delta_3'&=&\lim_{u\rightarrow 0} \frac{\mbox{d}}{\mbox{d}u}
\left\{\frac{u(3-p^{-2-2u})\zeta_q(1+2u)\zeta^3(2+2u)}{(1+p^{-2-2u})\zeta(4+4u)}
\left(\frac{4\pi^2 p^2}{q}\right)^{-u}\Gamma(1+u)^2e^{u^2}\right\}\nonumber\\
&=&\frac{(1-q^{-1})(3-p^{-2})}{1+p^{-2}}
\frac{\zeta^3(2)}{2\zeta(4)}
\left(\log{\frac{q}{4\pi^2p^2}}
+\frac{2\log p}{p^2+1}+\frac{2\log p}{3p^2-1}
+6\frac{\zeta'(2)}{\zeta(2)}-4\frac{\zeta'(4)}{\zeta(4)}\right)
+O\left(q^{-1}\log q\right).\nonumber\\
\eea
Thus by \eqref{D32}-\eqref{D33}, we obtain
\bna
\Delta_3=
\frac{(q-1)(3-p^{-2})}{24p(1+p^{-2})}
\frac{\zeta^3(2)}{\zeta(4)}
\left(\log{\frac{q}{4\pi^2p^2}}
+\frac{2\log p}{p^2+1}+\frac{2\log p}{3p^2-1}
+6\frac{\zeta'(2)}{\zeta(2)}-4\frac{\zeta'(4)}{\zeta(4)}\right)
+O\left(q^{\frac{1}{2}+\epsilon}\right).
\ena
Moreover, by \eqref{B3} and \eqref{E3}, we have
\bea \label{B30}
\mathcal
{B}_3(p^2,q,0)=
\frac{(q-1)(3-p^{-2})}{24p(1+p^{-2})}
\frac{\zeta^3(2)}{\zeta(4)}
\left(\log{\frac{q}{4\pi^2p^2}}
+\frac{2\log p}{p^2+1}+\frac{2\log p}{3p^2-1}
+6\frac{\zeta'(2)}{\zeta(2)}-4\frac{\zeta'(4)}{\zeta(4)}\right)
+O\left(q^{\frac{31}{32}+\epsilon}\right).\nonumber\\
\eea

If $t\neq 0$, there are simple poles at $u=0$ and $u=-it$. Then we have
\bea \label{D34}
\Delta_3'&=&\lim_{u\rightarrow -it}
\frac{(u+it)(3-p^{-2-2it-2u})\zeta_q(1+2it+2u)\zeta^3(2+2it+2u)}
{(1+p^{-2-2it-2u})\zeta(4+4it+4u)}
\left(\frac{4\pi^2 p^2}{q}\right)^{-u}\frac{\Gamma(1+it+u)^2}{u}e^{u^2}\nonumber\\
&+&\lim_{u\rightarrow 0}
\frac{(3-p^{-2-2it-2u})\zeta_q(1+2it+2u)\zeta^3(2+2it+2u)}
{(1+p^{-2-2it-2u})\zeta(4+4it+4u)}
\left(\frac{4\pi^2 p^2}{q}\right)^{-u}\Gamma(1+it+u)^2e^{u^2}\nonumber\\
&=&\frac{(1-q^{-1})(3-p^{-2})\zeta^3(2)}{2(1+p^{-2})\zeta(4)}\left(\frac{4\pi^2p^2}{q}\right)^{it}\frac{-e^{-t^2}}{it}
+\frac{(3-p^{-2-2it})\zeta_q(1+2it)\zeta^3(2+2it)}{(1+p^{-2-2it})\zeta(4+4it)}\Gamma(1+it)^2.\nonumber\\
\eea
Thus by \eqref{D32}, \eqref{E31} and \eqref{D34},
\bna
\Delta_3&=&\frac{q}{12p^{1+2it}}
\left(\frac{(1-q^{-1})(3-p^{-2})\zeta^3(2)}{2(1+p^{-2})\zeta(4)}\left(\frac{4\pi^2p^2}{q}\right)^{it}\frac{-e^{-t^2}}{it}\right.\\
&+&\left.\frac{(3-p^{-2-2it})\zeta_q(1+2it)\zeta^3(2+2it)}{(1+p^{-2-2it})\zeta(4+4it)}\Gamma(1+it)^2\right)
+O\left(q^{\frac{1}{2}+\epsilon}\right).
\ena
Moreover, by \eqref{B3} and \eqref{E3}, we have
\bea \label{B3p2qt}
\mathcal
{B}_3(p^2,q,t)&=&
\frac{q}{12p^{1+2it}}
\left(\frac{(1-q^{-1})(3-p^{-2})\zeta^3(2)}{2(1+p^{-2})\zeta(4)}\left(\frac{4\pi^2p^2}{q}\right)^{it}\frac{-e^{-t^2}}{it}\right.\nonumber\\
&+&\left.\frac{(3-p^{-2-2it})\zeta_q(1+2it)\zeta^3(2+2it)}{(1+p^{-2-2it})\zeta(4+4it)}\Gamma(1+it)^2\right)
+O\left(q^{\frac{31}{32}+\epsilon}\right).
\eea

\medskip

{\bf Evaluation of $\mathcal {B}_2(p^2,q,t)$}.
\medskip
We have
\bna
\mathcal {B}_2(p^2,q,t)&=&\sum_{d\geq 1 \atop (d,q)=1}\frac{1}{d^{1+2it}}
\sum_{p|n}\frac{\tau(n)}{n^{\frac{1}{2}+it}}W_t\left(\frac{4\pi^2 nd^2}{q}\right)
T_r(T_n).
\ena
By Lemma 2.2 and 2.3,
\bea \label{B21}
\mathcal {B}_2(p^2,q,t)&=&\frac{q}{12}\sum_{d\geq 1 \atop (d,q)=1}\frac{1}{d^{1+2it}}
\sum_{n =\diamond \atop p\mid n}\frac{\tau(n)}{n^{1+it}}W_t\left(\frac{4\pi^2 nd^2}{q}\right)\nonumber\\
&+&O\left(\sum_{B}
\left(B^{\frac{7}{4}}q^{-\frac{1}{2}}+B^{\frac{39}{32}+\epsilon}q^{\frac{1}{4}}\right)
\sum_{d\geq 1 \atop (d,q)=1}\frac{1}{d}
\max_{B\leq n\leq eB \atop p \mid n}\frac{\tau(n)}{\sqrt{n}}\left|
W_t\left(\frac{4\pi^2 nd^2}{q}\right)\right|\right)\nonumber\\
&:=&\Delta_2+E_2,
\eea
where the summation over $B$ is $\sum_{B=e^{\ell}, 0\leq \ell \leq [\log q^{11/10}]}$. By Lemma 2.2,
\bea \label{E2}
E_2
&\ll_{\epsilon}&q^{\epsilon}\sum_{ B}\left(B^{\frac{5}{4}}q^{-\frac{1}{2}}+B^{\frac{23}{32}+\epsilon}q^{\frac{1}{4}}\right)
\exp\left(-C\sqrt{\frac{4\pi^2 B}{q}}\right)\nonumber\\
&\ll_{\epsilon}&q^{\frac{31}{32}+\epsilon}.
\eea
By (2.3), we have
\bea \label{D2}
\Delta_2&=&\frac{q}{12}\sum_{d\geq 1 \atop (d,q)=1}\frac{1}{d^{1+2it}}
\sum_{n=\diamond \atop p\mid n}\frac{\tau (n)}{n^{1+it}}\frac{1}{2\pi i}\int\limits_{(2)}
\left(\frac{4\pi^2 nd^2}{q}\right)^{-u}\Gamma(1+it+u)^2e^{u^2}
\frac{\mbox{d}u}{u}\nonumber\\
&=&\frac{q}{12p^{2+2it}}\frac{1}{2\pi i}\int\limits_{(2)}\zeta_q(1+2it+2u)
\left(\sum_{\ell \geq1}\frac{\tau(p^2\ell^2)}{\ell^{2+2it+2u}}\right)
\left(\frac{4\pi^2p^2}{q}\right)^{-u}
\Gamma(1+it+u)^2e^{u^2}\frac{\mbox{d}u}{u}\nonumber\\
&=&\frac{1}{p}\Delta_3,
\eea
recalling the second display in \eqref{D31}.
By \eqref{B21}-\eqref{D2},
\bna
\mathcal {B}_2(p^2,q,t)=
\frac{1}{p}\Delta_3+O\left(q^{\frac{31}{32}+\epsilon}\right).
\ena
Therefore, if $ t=0 $,
\bea \label{B20}
\mathcal
{B}_2(p^2,q,0)
=\frac{(q-1)(3-p^{-2})}{24p^2(1+p^{-2})}
\frac{\zeta^3(2)}{\zeta(4)}
\left(\log{\frac{q}{4\pi^2p^2}}
+\frac{2\log p}{p^2+1}+\frac{2\log p}{3p^2-1}
+6\frac{\zeta'(2)}{\zeta(2)}-4\frac{\zeta'(4)}{\zeta(4)}\right)
+O\left(q^{\frac{31}{32}+\epsilon}\right).\nonumber\\
\eea
If $t\neq 0$,
\bea \label{B2p2qt}
\mathcal
{B}_2(p^2,q,t)&=&
\frac{q}{12p^{2+2it}}
\left(\frac{(1-q^{-1})(3-p^{-2})\zeta^3(2)}{2(1+p^{-2})\zeta(4)}\left(\frac{4\pi^2p^2}{q}\right)^{it}\frac{-e^{-t^2}}{it}\right.\nonumber\\
&+&\left.\frac{(3-p^{-2-2it})\zeta_q(1+2it)\zeta^3(2+2it)}{(1+p^{-2-2it})\zeta(4+4it)}\Gamma(1+it)^2\right)
+O\left(q^{\frac{31}{32}+\epsilon}\right).
\eea

\medskip

{\bf Conclusion}. By \eqref{BB}, \eqref{B10}, \eqref{B30} and \eqref{B20}, we have for $t=0$,
\bna
\mathcal
{B}(p^2,q,0)
&=&\frac{(q-1)\zeta^3(2)}{24p^2\zeta(4)}
\left(p+\frac{(p+1)(3-p^{-2})}{1+p^{-2}}\right)\left(\log\frac{q}{4\pi^2p^2}
+6\frac{\zeta'(2)}{\zeta(2)}
-4\frac{\zeta'(4)}{\zeta(4)}\right)\\
&+&\frac{(q-1)\zeta^3(2)\log p}{12p\zeta(4)}
\left(1+\frac{4p(p+1)(3-p^{-2})}{(1+p^{-2})(p^2+1)(3p^2-1)}
\right)\\
&+&O\left(p^{\frac{39}{16}+\epsilon}q^{\frac{31}{32}+\epsilon}\right),
\ena
and by \eqref{Ap2qt},
\bna
\mathcal {A}(p^2,q,0)&=&
\frac{(q-1)\zeta^3(2)}{12p^2\zeta(4)}
\left(p+\frac{(p+1)(3-p^{-2})}{1+p^{-2}}\right)\left(\log\frac{q}{4\pi^2p^2}
+6\frac{\zeta'(2)}{\zeta(2)}
-4\frac{\zeta'(4)}{\zeta(4)}\right)\\
&+&\frac{(q-1)\zeta^3(2)\log p}{6p\zeta(4)}
\left(1+\frac{4p(p+1)(3-p^{-2})}{(1+p^{-2})(p^2+1)(3p^2-1)}
\right)\\
&+&O_\epsilon\left(p^{\frac{39}{16}+\epsilon}q^{\frac{31}{32}+\epsilon}\right),
\ena
which is an asymptotic formula for $p<q^{\frac{1}{110}-\epsilon}$.

Let $0\neq t\in \mathbb{R}$ be fixed.
By \eqref{BB}, \eqref{B1p2qt}, \eqref{B3p2qt} and \eqref{B2p2qt}, we have
\bna
\mathcal
{B}(p^2,q,t)&=&
\frac {q}{24p^{2+2it}}
\frac{(1-q^{-1})\zeta^3(2)}{\zeta(4)}
\left(\frac{4\pi^2p^2}{q}\right)^{it}
\left(p+\frac{(p+1)(3-p^{-2})}{1+p^{-2}}\right)
\frac{-e^{-t^2}}{it}\nonumber\\
&+&
\frac {q}{12p^{2+2it}}
\frac{\zeta_q(1+2it)\zeta^3(2+2it)}{\zeta(4+4it)}
\left(p^{1+2it}+\frac{(p+1)(3-p^{-2-2it})}{1+p^{-2-2it}}\right)
\Gamma(1+it)^2\\
&+&O(p^{\frac{39}{16}+\epsilon}q^{\frac{31}{32}+\epsilon}),
\ena
and by \eqref{Ap2qt},
\bna
\mathcal {A}(p^2,q,t)&=&
\frac{\zeta_q(1+2it)\zeta^3(2+2it)}{12\zeta(4+4it)}
\left(\frac{q}{p}+\frac{q(p+1)(3-p^{-2-2it})}{p^{2+2it}+1}\right)\\
&+&\frac{\Gamma(1-it)^2}{\Gamma(1+it)^2}
\left(\frac{q}{4\pi^2}\right)^{-2it}
\frac{\zeta_q(1-2it)\zeta^3(2-2it)}{12\zeta(4-4it)}
\left(\frac{q}{p}+\frac{q(p+1)(3-p^{-2+2it})}{p^{2-2it}+1}\right)\\
&+&O_{t,\epsilon}\left(p^{\frac{39}{16}+\epsilon}q^{\frac{31}{32}+\epsilon}\right),
\ena
which is an asymptotic formula for $p<q^{\frac{1}{110}-\epsilon}$.

\vskip0.2in

\bigskip

{
Wei Liu

School of Mathematics and Statistics

Shandong University at Weihai

Weihai, Shandong 264209

China

{\it Email: wei.liu@mail.sdu.edu.cn}}

\end{document}